\documentclass[12pt]{article}
\usepackage{amssymb,latexsym}
\def\QQ{{\mathbb Q}}
\def\RR{{\mathbb R}}

\newtheorem{formula}{}[section]
\newtheorem{proposition}[formula]{Proposition}
\newtheorem{definition}[formula]{\indent Definition}
\newtheorem{corollary}[formula]{\indent Corollary}
\newtheorem{remark}[formula]{\indent Remark}
\newtheorem{lemma}[formula]{\indent Lemma}
\newtheorem{theorem}[formula]{\indent Theorem}

\newtheorem{example}[formula]{Example}

\def\thrm{\begin{theorem}}
\def\thrml#1{\begin{theorem}\label{#1}}
\def\ethrm{\end{theorem}}
\def\rmrk{\begin{remark}}
\def\rmrkl#1{\begin{remark}\label{#1}}
\def\ermrk{\end{remark}}
\def\dfntn{\begin{definition}}
\def\dfntnl#1{\begin{definition}\label{#1}}
\def\edfntn{\end{definition}}
\def\nmrt{\begin{enumerate}}
\def\enmrt{\end{enumerate}}

\def\qtn{\begin{equation}}
\def\qtnl#1{\begin{equation}\label{#1}}
\def\eqtn{\end{equation}}
\def\lmm{\begin{lemma}}
\def\lmml#1{\begin{lemma}\label{#1}}
\def\elmm{\end{lemma}}
\def\crllr{\begin{corollary}}
\def\crllrl#1{\begin{corollary}\label{#1}}
\def\ecrllr{\end{corollary}}

\begin{document}
\title{}
\date{}
\maketitle
%\nopagenumbers
%\begin{titlepage}
%\title
\vspace{-0,1cm} \centerline{\bf Decomposing tropical rational functions}
%\centerline{\bf AND RELATIVE KOLMOGOROV COMPLEXITY}
\vspace{7mm}
\author{
\centerline{Dima Grigoriev}
%\\[-1pt]
\vspace{3mm}
%\centerline{$^1$ St.Petersbourg University,
% Universitetskaya nab., 7/9,}
%\centerline{ St.Petersbourg,
% 199164, RUSSIA}
%\vspace{3mm}
%Dima Grigoriev \\[-1pt]
\centerline{CNRS, Math\'ematique, Universit\'e de Lille, Villeneuve
d'Ascq, 59655, France} \vspace{1mm} \centerline{e-mail:\
dmitry.grigoryev@univ-lille.fr } \vspace{1mm}
\centerline{URL:\ http://en.wikipedia.org/wiki/Dima\_Grigoriev} }
%\date{}
%\maketitle

\begin{abstract}
An algorithm is designed which decomposes a tropical univariate rational function into a composition
of tropical binomials and trinomials. When a function is monotone, the composition consists
just of binomials. Similar algorithms are designed for decomposing tropical algebraic rational
functions being (in the classical language) piece-wise linear functions with rational slopes of
their linear pieces. In addition, we provide a criterion when the composition of two tropical polynomials commutes (for classical polynomials a similar question was answered by J.~Ritt).
\end{abstract}

\section*{Introduction}

We study decomposing tropical univariate rational functions (compositions of tropical rational functions find applications in deep learning of neural networks, see e.~g. \cite{neural}). A tropical rational function is the tropical
quotient (which corresponds to the subtraction in the classical sense) of two tropical polynomials.
Thus, a tropical rational function is (classically) a piece-wise linear function with integer slopes
of its linear pieces. A tropical root of a tropical rational function is defined as a point at which
the function is not differentiable.

Relaxing the requirement that the slopes are integers allowing them to be rationals, we arrive to the
concept of {\it tropical algebraic rational functions} or tropical Newton-Puiseux rational functions \cite{G}
playing the role of algebraic functions in tropical algebra.

In classical algebra the problem of decomposing polynomials, rational and algebraic functions was
elaborated in \cite{Gathen}, \cite{Kozen}, \cite{Rubio}. In tropical algebra the answer to the decomposing problem differs essentially from its classical counterpart. We show that a tropical  rational function is a
composition of binomials and trinomials. The similar holds for tropical algebraic rational functions.

In section~\ref{one} we introduce tropical monotone rational and algebraic rational functions and bound the
number of tropical roots of compositions of tropical polynomials, monotone rational functions and rational functions.

In section~\ref{two} we design an algorithm which decomposes a tropical algebraic function and also a tropical monotone algebraic rational function into a composition of tropical binomials. In addition, we design an algorithm which decomposes a tropical algebraic rational function into a composition of tropical binomials and trinomials. Moreover, we provide a bound on the number of composants.

In section~\ref{three} decompositions of tropical rational functions  (so, with integer slopes of their linear pieces) are studied. We design an algorithm which decomposes a tropical monotone rational function into a composition of tropical monotone binomials and monotone trinomials. Also we design an algorithm which decomposes a tropical rational function into a composition of tropical binomials and trinomials. In addition, a criterion is provided, when a tropical monotone trinomial is decomposable. Finally, bounds on the number of composants are given. 
 
In section~\ref{four} we prove that the composition of two tropical polynomials $f,\, g$ without free terms commutes: $f\circ g=g\circ f$ iff there is a common fixed point $x_0$ (perhaps, $x_0=\infty$) for both $f,\, g$, i.~e. $f(x_0)=g(x_0)=x_0$ and there exist a tropical increasing algebraic rational function $h$ and integers $a,\, b\ge 1,\, k,\, m\ge 0$ such that either $f=h^k,\, g=h^m$ or $f=ax+x_0(1-a),\, g=bx+x_0(1-b)$  
%(respectively, $f=h^p,\, g=h^q$) 
on the interval $(-\infty, x_0]$ (similar conditions hold on the interval $[x_0,\, \infty)$), unless $f=x+c_1,\, g=x+c_2,\, x\in \RR$ for some $c_1,\, c_2 \in \RR$. Also we provide an example of a pair of 
increasing tropical rational functions with commuting composition which do not satisfy the latter conditions.
For classical polynomials the answer to commutativity was given in \cite{R22}, \cite{R23} (more recent generalizations and further references one can find in \cite{P}), in which commuting Chebyshev polynomials play a crucial role.

In section~\ref{five} we introduce tropical polynomial (respectively, Laurent polynomial and rational) parametrizations of polygonal lines. We show that any polygonal line admits a tropical rational parametrization and provide criteria when it does admit a tropical polynomial (respectively, Laurent polynomial) parametrization.

\section{Tropical monotone rational functions}\label{one}

Recall (see e.~g. \cite{MS}) that a (univariate) tropical polynomial has a
form $f=\min_{0\le i\le d}\{ix+a_i\},\, a_i\in \RR \cup \{\infty\}$. Linear
functions $ix+a_i,\, 0\le i\le d$ are called tropical monomials. So, the 
minimum plays the role of the addition in tropical algebra, while the addition
plays the role of the multiplication. Thus, $f$ is a convex piece-wise linear
function with integer slopes of the edges of its graph (sometimes, slightly 
abusing the terminology we call them the edges of $f$). We consider the natural ordering of the edges 
from the left to the right. A point $x\in \RR$ is a
tropical root of $f$ if the minimum in $f$ is attained at least for two linear functions
$ix+a_i,\, 0\le i\le d$. In other words, tropical roots of $f$ are the points at which $f$ is not
differentiable.

A tropical rational function is a difference (which plays the role of the division in tropical algebra) of
two tropical polynomials. It is a piece-wise linear function. So, its graph consists of several edges.
Conversely, any continuous piece-wise linear function with integer slopes of its linear pieces (edges of its graph) is a tropical rational function (cf. \cite{G} where one can find further references). As the roots of a tropical rational function we again mean the points at which the function is not differentiable.

If $g,\, h$ are tropical rational functions with $p,\, q$ tropical roots, respectively, then (see \cite{GP})
the number of the roots of 

$\bullet$ $\min\{g,\, h\}$ is at most $p+q+1$;

$\bullet$ $g+h$ or $g-h$ is at most $p+q$. 

In this paper we study compositions $g\circ h$ (being tropical rational functions as well). Note that if $g,\, h$ are tropical polynomials then  $g\circ h$ is also a tropical polynomial. If $s_1,\dots,s_k$ are consecutive (integer) slopes of (the linear pieces of) a tropical rational function $g$ then $g$ is a tropical polynomial iff $s_1>\cdots >s_k\ge 0$. 

In a tropical {\it monotone} increasing (or decreasing, respectively)  rational function $g$ its slopes are positive (respectively, negative). Note that any tropical polynomial $\min_{1\le i\le d} \{ix+a_i\}$ without free term is monotone increasing. 

One can directly verify the following proposition.

\begin{proposition}\label{roots}
If $g,\, h$ are tropical monotone  rational functions with $p,\, q$ tropical roots, respectively,  then the tropical monotone rational function $g\circ h$ has at most $p+q$ tropical roots. Moreover, if on an interval
$[a,b]\subset \RR$ function $h$ is linear with a slope $s$, and $g$ is linear with a slope $l$ on the interval $[h(a),\, h(b)]$ (respectively, $[h(b),\, h(a)]$) when $h$ increases (respectively, decreases) then
$g\circ h$ is linear on the interval $[a,\, b]$ with the slope $sl$.
\end{proposition} 

\begin{remark}\label{roots-rational}
In general, the number of tropical roots of the composition $g\circ h$ of tropical rational functions does not exceed $pq+p+q$. Moreover, if $s_0,\dots,s_p$ (respectively, $t_0,\dots,t_q$) are the slopes of (the graph of) $g$ (respectively, $h$) listed with possible repetitions (multiplicities), then the slopes of $g\circ h$ are among $s_it_j,\, 0\le i\le p,\, 0\le j\le q$.

For a tropical rational function $g=\min\{x+1,\, -x+1\}$ the number of the tropical roots of $k$ iterations of $g^k:=g\circ\cdots\circ g$ is $2^k-1$ \cite{neural} (see also \cite{GP}).
\end{remark}

Admitting rational coefficients in $\min_i \{b_ix+a_i\},\, 0\le b_i\in \QQ$, we arrive to the concept of {\it tropical algebraic functions} (or {\it tropical Newton-Puiseux polynomials}) \cite{G}. Respectively, we consider {\it tropical algebraic rational functions} being differences of tropical algebraic functions \cite{G}. 

\begin{remark}
The above statement in Proposition~\ref{roots} on the slopes of tropical rational functions holds for tropical algebraic rational functions as well with the difference that now we admit rational slopes rather than just integers. The above bounds on the number of tropical roots also hold literally for tropical algebraic rational functions.
\end{remark}

\section{Decomposing tropical algebraic rational functions}\label{two}

In this section we consider tropical algebraic rational functions. As a {\it tropical algebraic rational binomial} we mean a function of the form either
$\min\{b_1x+a_1,\, b_2x+a_2\},\, 0\neq b_1,\, b_2 \in \QQ$ or $\max\{b_1x+a_1,\, b_2x+a_2\}$. In the geometric language the former function is a convex piece-wise linear function with two (unbounded) edges (and we call it a {\it tropical algebraic binomial}), while the latter one is concave.
If $b_1,\, b_2 >0$ then in both cases the functions are monotone increasing.

\begin{proposition}\label{monotone-algebraic}
(i) There is an algorithm which for a tropical algebraic function $f$ with $k$ tropical roots yields a decomposition of $f$ into $k$ tropical  algebraic binomials;

(ii) let $f$ be a tropical monotone algebraic rational function with $k$ tropical roots. Then the algorithm yields a decomposition of $f$ into $k$ tropical monotone algebraic rational binomials.
\end{proposition}

\begin{remark}\label{monotone-algebraic-number}
Due to Proposition~\ref{roots} and taking into the account that each tropical algebraic rational binomial has a single tropical root, we conclude that in Proposition~\ref{monotone-algebraic} one can't take less than $k$ composants.
\end{remark}

{\bf Proof}. The proofs for both items (i), (ii) proceed similarly. Let $f$ have consecutive slopes $s_0,\dots,s_k$ of its linear pieces. Recall that $s_0>s_1>\cdots >s_k\ge 0$ in case   (i) and $s_1,\dots,s_k>0$ in case  (ii). Denote by $x_l$ the $l$-th tropical root of $f,\, 1\le l\le k$. Take a (piece-wise linear) function $h$ with $k$ slopes
$$s_0.s_1,\dots,s_{l-1},s_{l+1}\cdot s_{l-1}/s_l,\dots,s_k\cdot s_{l-1}/s_l$$
\noindent coinciding with $f$ for $x\le x_l$ and replacing $f$ by the composition with the linear function $((s_{l-1}/s_l)x+f(x_l)(1-s_{l-1}/s_l))\circ f$ for $x\ge x_l$. Thus, $h$ has the tropical roots $x_1,\dots,x_{l-1},x_{l+1},\dots,x_k$. The described procedure replacing $f$ by $h$ we call {\it straightening}:  one tropical root (at $x_l$) disappears. 

Take a tropical algebraic rational binomial $g$ coinciding with the identity function $x\to x$ for $x\le f(x_l)$ and with the linear function $(s_l/s_{l-1})x+f(x_l)(1-s_l/s_{l-1})$ for $x\ge f(x_l)$. Then $f=g\circ h$. Note that in case (i) $g$ is a tropical algebraic binomial since $s_{l-1}/s_l<1$.

Proceeding by induction on $k$ we complete the proof of the Proposition. $\Box$

\begin{remark}
Observe that each tropical root of $f$ corresponds to a suitable composant in a decomposition of $f$. Thus, by choosing (in the proof of Proposition~\ref{monotone-algebraic} above) the tropical roots in different orders, we obtain $k!$ "combinatorially different types" of decompositions of $f$.
\end{remark}

Now let $f$ be a tropical algebraic rational function. Let for definiteness the first edge of $f$ with a non-zero slope have a positive slope. Consider tropical roots $x$ of $f$ such that $f$ has an edge with a negative slope to the right of $x$. If there does not exist such $x$ then $f$ is (non-strictly) monotone increasing, and we proceed to study the monotone case later. Among such $x$ pick $x_0$ (perhaps, if not unique then pick any of them) with the maximal value $f(x_0)$. Take a tropical root $x_1>x_0$ of $f$ with the minimal value $f(x_1)$. We have $f(x_1)<f(x_0),\, f([x_0,\, x_1])=[f(x_1),\, f(x_0)]$ and $\max\{f(x):\, x\in (-\infty,\, x_0]\}=f(x_0)$.  

First we consider the case when $f(x)\le f(x_0)$ for all $x\ge x_1$. Note that in this case $f(x)\le f(x_0)$ for all $x\ge x_0$ due to the choice of $x_0$. 

If both adjacent to $x_0$ edges of $f$ have non-zero slopes then the edge to the left from $x_0$  has a positive slope $s_0>0$, while the edge to the right from $x_0$ has a negative slope $s_1<0$ (due to the choice of $x_0$). Take as $g$ a tropical algebraic rational binomial which coincides with the identity function $x\to x$ for $x\le f(x_0)$ and with a linear function $(s_1/s_0)x+f(x_0)(1-s_1/s_0)$ for $x\ge f(x_0)$. So, $g$ is a tropical non-monotone algebraic rational binomial.

As $h$ take a tropical algebraic rational function which coincides with $f$ for $x\le x_0$ and coincides with the composition with the linear function $((s_0/s_1)x+f(x_0)(1-s_0/s_1))\circ f$ for $x\ge x_0$. Then $f=g\circ h$. By a {\it block of edges} of $f$ we mean a sequence of consecutive edges of the equal signs of their slopes (ignoring edges with zero slopes). Observe that $h$ has one less block of edges than $f$ does. Thus, by passing from $f$ to $h$ we straighten $f$ at point $x_0$.

Otherwise, if one of adjacent to $x_0$ edges has zero slope then as $g$ take a tropical binomial coinciding with the identity function $x\to x$ for $x\le f(x_0)$ and with the linear function $-x+2f(x_0)$ for $x\ge f(x_0)$. As $h$ take a tropical algebraic rational function which coincidies with $f$ for $x\le x_0$ and with the composition with the linear function $(-x+2f(x_0))\circ f$ for $x\ge x_0$. Then again $f=g\circ h$, and $h$ has one less block of edges than $f$ does.
On the other hand, $h$ has the same number of tropical roots as $f$ does, so one does not straighten a piece-wise linear function at a point if one of two adjacent edges to this point has zero slope.

Now we proceed to the case when $f(x)$ takes a value greater than $f(x_0)$ for some $x>x_1$. Then $\min\{f(x): \, x\ge x_0\} = f(x_1)$. 

A {\it tropical regular algebraic rational trinomial} is a piece-wise linear function with 3 edges having rational non-zero slopes. If the slopes are decreasing or increasing positive integers we talk about a {\it tropical trinomial}. 

Construct the following tropical algebraic rational functions $h,\, g$. If both the edge of $f$ with the right (and respectively, the left) end-point $(x_0,f(x_0))$ has a non-zero slope $s_+$ (respectively, a non-zero slope $s_-$) then $h$ on the interval $(-\infty,\, x_0]$ coincides with the composition $(-(s_-/s_+)x+f(x_0)(1+s_-/s_+))\circ f$. Note that $s_+>0,\, s_-<0$. As $g$ take a function which on the interval $(-\infty,\, f(x_0)]$ coincides with the linear function $-(s_+/s_-)x+f(x_0)(1+s_+/s_-)$. We have $\max\{h(x): \, x\le x_0\}=f(x_0)$ and $g((-\infty,\, f(x_0)])=(-\infty,\, f(x_0)]$. In case if $s_+\cdot s_- =0$ then $h$ on the interval $(-\infty,\, x_0]$ coincides with $f$, and $g$ on the interval $(-\infty,\, f(x_0)]$ coincides with the identity function $x\to x$.

On the interval $[x_0,\, x_1]$ the function $h$ in both cases coincides with the composition
$(-x+2f(x_0))\circ f$, and $g$ on the interval $[f(x_0),\, 2f(x_0)-f(x_1)]$ coincides with the linear function $-x+2f(x_0)$. Then $h([x_0,\, x_1])=[f(x_0),\, 2f(x_0)-f(x_1)]$ and $g([f(x_0,\, 2f(x_0)-f(x_1)])=[f(x_1),\, f(x_0)]$.  

Finally, define $h$ on the interval $[x_1,\infty)$ and $g$ on the interval $[2f(x_0)-f(x_1),\, \infty)$. Similar to the consideration above of the interval $(-\infty, x_0]$ denote by $t_-$ (respectively, $t_+$) the slope of the edge of $f$ with the right (respectively, the left) end-point $(x_1,\, f(x_1))$. If $t_-\cdot t_+ \neq 0$ (in this case $t_-<0,\, t_+>0$ due to the choice of $x_1$) then $h$ on the interval $[x_1,\infty)$ coincides with the composition with the linear function $(-(t_-/t_+)x+2f(x_0)+f(x_1)(t_-/t_+-1))\circ f$. In this case $g$ on the interval $[2f(x_0)-f(x_1),\, \infty)$ coincides with the linear function $-t_+/t_-x + f(x_1)+t_+/t_-(2f(x_0)-f(x_1))$. Then $\min\{h(x)\, : \, x_1\le x<\infty\}=
2f(x_0)-f(x_1)$ and $g([2f(x_0)-f(x_1),\, \infty))=[f(x_1),\, \infty)$.

Otherwise, if $t_-\cdot t_+ =0$ then $h$ on the interval $[x_1,\, \infty)$ coincides with the composition $(x+2f(x_0)-2f(x_1))\circ f$, and $g$ on the interval  $[2f(x_0)-f(x_1),\, \infty)$ coincides with the linear function $x-2f(x_0)+2f(x_1)$. In this case again $\min\{h(x)\, :\, x_1\le x<\infty\}=2f(x_0)-f(x_1)$ and
$g([2f(x_0)-f(x_1),\, \infty))=[f(x_1),\, \infty)$. Thus, $f=g\circ h$.

Observe that $h$ has two less blocks of edges than $f$ does. Also note that $x$ is a tropical root of $h$ with the adjacent to $x$ edges of $h$ with the equal signs of their slopes iff $x$ is a tropical root of $f$ satisfying the same property (we call such $x$ a non-extremal tropical root of $h$ because $x$ is not a local extremal of $h$). In addition, the numbers of edges with zero slope are the same for $f$ and for $h$.

Thus, applying two described decomposition procedures to $f$ and obtaining $g$ to be either a tropical non-monotone algebraic rational binomial or a tropical regular algebraic rational trinomial, while it is possible, we arrive to a tropical algebraic rational function $f_0$ which is non-decreasing, so the slopes of its edges are non-negative. Thus, 
\begin{eqnarray}\label{1}
f=g_1\circ \cdots \circ g_k \circ f_0
\end{eqnarray} 
\noindent where each of $g_1,\dots,g_k$ is either a tropical non-monotone algebraic rational binomial (their number among $g_1,\dots,g_k$ denote by $k_2$) or a tropical regular algebraic rational trinomial (their number denote by $k_3:=k-k_2$). Therefore, $k_2+2k_3$ equals the number of blocks of edges of $f$.

Now take a non-extremal tropical root $x_4$ of $f_0$. Let $s_->0$ (respectively, $s_+>0$) be the slope of the adjacent to $x_4$ left edge (respectively, right edge) of $f_0$. Denote by $g^{(1)}$ a tropical monotone increasing algebraic rational binomial which coincides with the identity function on the interval $(-\infty,\, f_0(x_4)]$ and which coincides with the linear function $(s_+/s_-)x+f(x_0)(1-s_+/s_-)$ on the interval $[f_0(x_4),\, \infty)$. Denote by $h_0$ a tropical non-decreasing algebraic rational function which coincides with $f_0$ on the interval $(-\infty,\, x_4]$ and which coincides with the composition $((s_-/s_+)x+f_0(x_4)(1-s_-/s_+))\circ f_0$ on the interval $[x_4,\, \infty)$.

Then $f_0=g^{(1)}\circ h_0$ and $h_0$ has no tropical root at $x_4$, while having all other tropical roots of $f_0$, so $h_0$ is a straightening of $f_0$. Applying the just described procedure to all non-extremal tropical roots of $f_0$, we obtain a decomposition 
\begin{eqnarray}\label{2}
f_0=g^{(1)}_1\circ \cdots \circ g^{(1)}_{k_1}\circ f^{(1)}
\end{eqnarray}
\noindent where each of $g^{(1)}_1,\dots,g^{(1)}_{k_1}$ is a tropical increasing algebraic rational binomial. Every second edge of $f^{(1)}$ has zero slope, and the number $k_0$ of edges with zero slope of $f^{(1)}$ equals the same number of $f$. Observe that $k_1$ equals the number of non-extremal tropical roots of $f$ (and also equals the number of non-extremal tropical roots of $f_0$). Hence $k_1+k_2+2k_3$ does not exceed the number of edges with non-zero slopes of $f$.

As a {\it tropical singular algebraic rational trinomial} we mean a trinomial whose middle edge has zero slope. Slightly abusing the terminology, we admit singular trinomials without one or two edges with non-zero slopes.

We are looking for a decomposition 
\begin{eqnarray}\label{3}
f^{(1)}=g^{(0)}_{k_0}\circ \cdots \circ g^{(0)}_1
\end{eqnarray} 
\noindent where $g^{(0)}_i,\, 1\le i\le k_0$ is a tropical singular algebraic trinomial. To decompose take the left-most interval $[x_0,\, x_1]$ on which $f^{(1)}$ is constant, in other words, the edge of $f^{(1)}$ on $[x_0,\, x_1]$ has zero slope. It can happen that $x_0=-\infty$, in this case some of the following considerations become void. Define $g^{(0)}_1$ on the interval $(-\infty,\, x_0]$ as the identity function and on the interval $[x_0,\, x_1]$ as the constant function with the value $x_0$. 

Let $f^{(1)}$ on the interval  $(-\infty,\, x_0]$ equal a linear function $sx+r$ (so, $s$ is the slope of the left-most edge of $f^{(1)}$), in particular $sx_0+r=f^{(1)}(x_0)$. Define $g^{(0)}$ (later we'll get that $g^{(0)}=g^{(0)}_{k_0}\circ \cdots \circ g^{(0)}_2$) on the interval $(-\infty,\, x_0]$ as the linear function $sx+r$. Therefore, $g^{(0)}\circ g^{(0)}_1$ on the interval $(-\infty,\, x_0]$ coincides with $f^{(1)}$. The same coincidence holds on the interval $[x_0,\, x_1]$ as well. Let the edge of $f^{(1)}$ with the left end-point $(x_1,\, f^{(1)}(x_1)=sx_0+r)$ have a slope $p$. Then define $g^{(0)}_1$ on the interval $[x_1,\, \infty)$ as the linear function $(p/s)x+x_0-(p/s)x_1$. Also define $g^{(0)}$ on the interval $[x_0,\, \infty)$ as
the composition $f^{(1)}\circ ((s/p)x-(s/p)x_0+x_1)$. Then $f^{(1)}=g^{(0)}\circ g^{(0)}_1$.

Now we observe that $g^{(0)}$ is a continuous non-decreasing piece-wise linear function: we have constructed it by gluing at $x_0$ two non-decreasing piece-wise linear functions both having the value $sx_0+r=f^{(1)}(x_0)=f^{(1)}(x_1)$ at $x_0$. Moreover, the slope of the edge of $g^{(0)}$ with the right end-point $(x_0,\, f^{(1)}(x_0))$ equals $s$ which coincides with the slope of the edge of  $g^{(0)}$ with the left end-point $(x_0,\, f^{(1)}(x_0))$. Therefore, $g^{(0)}$ has no tropical root at $x_0$ (so, $g^{(0)}$ is a straightening of $f^{(1)}$), and $g^{(0)}$ is of a similar shape as $f^{(1)}$, i.~e. $g^{(0)}$ is a non-decreasing piece-wise linear function whose every second edge has zero slope. On the other hand,  $g^{(0)}$ has one less edge with zero slope than $f^{(1)}$ does. Continuing in this way, we construct a required decomposition (\ref{3}).

Combining (\ref{1}), (\ref{2}) and (\ref{3}) we complete the proof of the following theorem.

\begin{theorem}\label{algebraic-rational}
There is an algorithm which decomposes a tropical algebraic rational function 
\begin{eqnarray}\label{4}
f=g_1\circ \cdots \circ g_k\circ g^{(1)}_1\circ \cdots \circ g^{(1)}_{k_1} \circ g^{(0)}_{k_0} \circ \cdots \circ g^{(0)}_1
\end{eqnarray}
\noindent where each $g_i,\, 1\le i\le k$ is either a tropical regular algebraic rational trinomial or a tropical non-monotone algebraic rational binomial (cf. (\ref{1})), each $g^{(1)}_j,\, 1\le j\le k_1$ is a tropical monotone algebraic rational binomial (cf. (\ref{2})), and each $g^{(0)}_l,\, 1\le l\le k_0$ is a tropical singular algebraic trinomial (cf. (\ref{3})). 

Moreover, if $k_3$ is the number of tropical regular algebraic rational trinomials, and $k_2$ is the number of tropical non-monotone algebraic rational binomials in (\ref{4}), so $k_3+k_2=k$ then $2k_3+k_2$ is the number of blocks of edges of $f$ of the equal (non-zero) signs of their slopes. The number $2k_3+k_2+k_1$ does not exceed the number of edges of $f$ with non-zero slopes, finally $k_0$ equals the number of edges with zero slopes.
\end{theorem}

\begin{remark}
The number of tropical roots of $f$ is greater or equal to $k_1+2k_0$, and on the other hand, is less or equal to $2k_3+k_2+k_1+k_0$, the latter number also equals the total number of tropical roots in the composants of $f$ from (\ref{4}) (cf. Remark~\ref{roots-rational}).
\end{remark}

\section{Tropical rational functions}\label{three}

In this section we study decompositions of tropical rational functions,  we recall that the slopes of edges of a piece-wise rational function $f$  are integers (unlike the section~\ref{two} in which the slopes could be rationals).

\begin{theorem}\label{tropical-rational}
(i) There is an algorithm which for a tropical monotone rational function $f$ yields its decomposition into tropical monotone binomials and tropical monotone trinomials. The number of composants does not exceed the number of tropical roots of $f$ (cf. Proposition~\ref{monotone-algebraic} and Remark~\ref{monotone-algebraic-number});

(ii) there is an algorithm which decomposes a tropical rational function 
%\begin{eqnarray}\label{5}
$$f=g_1\circ \cdots \circ g_k\circ h_1\circ \cdots \circ h_m\circ g^{(0)}_{k_0}\circ \cdots \circ g^{(0)}_1$$
%\end{eqnarray}
\noindent (cf. (\ref{4})) where each $g_i,\, 1\le i\le k$ is either a tropical non-monotone  rational binomial with $\pm 1$ slopes or a tropical non-monotone rational trinomial with $\pm 1$ slopes (cf. (\ref{1})),  each $h_j,\, 1\le j\le m$ is either a tropical regular monotone binomial or a tropical regular monotone trinomial, and each $g^{(0)}_l,\, 1\le l\le k_0$ is a tropical singular monotone trinomial. The number of binomials among $g_1,\dots,g_k$ plus the double number of trinomials among $g_1,\dots,g_k$ does not exceed the number of blocks of edges of $f$  (cf. Theorem~\ref{algebraic-rational}). The number $m$ does not exceed the number of edges of $f$, and the number $k_0$ equals the number of edges of $f$ with zero slopes (again cf. Theorem~\ref{algebraic-rational} and (\ref{3}));

(iii) let $f$ be a tropical monotone rational function (respectively, a tropical polynomial) with the slopes of its edges $a_0,\dots,a_n\ge 1$ (respectively, $a_0>\dots >a_n\ge 1$) and denote by $q_i,\, 1\le i\le n$ the denominator of the irreducible fraction $a_i/a_{i-1}$. Then $f$ is a composition of tropical rational binomials (respectively, tropical binomials) iff $(q_1\cdots q_n)|a_0$.

%(iii) a tropical regular monotone trinomial with (integer) slopes $a,\, b,\, c >0$ is decomposable (in%to a pair of tropical regular monotone binomials) iff $b|(a\cdot c)$. 
\end{theorem}

\begin{remark}
If $f$ satisfies the latter condition in (iii) we call $f$ {\it completely decomposable}. 

This condition in (iii) is equivalent to a more symmetric one: for any $m\ge 1$ and $j\ge 0$ such that $m+2j<n$ it holds
$$\prod_{0\le i\le j} a_{m+2i} \quad | \prod_{0\le i\le j+1} a_{m+2i-1}.$$
\noindent In particular, for $n=2$ (trinomials), the condition in (iii) for $a_0,\, a_1,\, a_2$ is equivalent to $a_1|(a_0a_2)$.
%In case of tropical monotone rational trinomials (as well as in case of tropical trinomials) with slopes $a_0,\, a_1,\, a_2$ the condition in (iii) is equivalent to $a_1|(a_0a_2)$. 
\end{remark}

{\bf Proof}. (i). If an increasing $f$ has at least 4 edges then take any its tropical root $x_0$ being neither the left-most nor the right-most. Define $h$ to coincide with $f$ on the interval $(-\infty,\, x_0]$ and $g$ to coincide with the identity function on the interval $(-\infty,\, f(x_0)]$. Then define $h$ on the interval $[x_0,\, \infty)$ to coincide with the linear function $x+f(x_0)-x_0$, and define $g$ on the interval $[f(x_0),\, \infty)$ to coincide with the composition $f\circ (x-f(x_0)+x_0)$. Then $f=g\circ h$.

Continuing in this way, applying further the described construction to $g,\, h$ we complete the proof of (i).

(ii). First, similar to the proof of Theorem~\ref{algebraic-rational} one represents (by means of straightening) $f=g\circ h$ (assume w.l.o.g. that the first edge of $f$ with non-zero slope has a positive slope), where $g$ is either a tropical non-monotone binomial with the slopes of its edges $1$ and $-1$ or a tropical trinomial with the slopes $1,\, -1,\, 1$, while $h$ being a tropical rational function with less number of blocks of edges than $f$.

Continuing in this way, while it is possible, we arrive to a tropical non-decreasing rational function $f^{(0)}$ such that $f=g_1\circ \cdots \circ g_k \circ f^{(0)}$ (cf. (\ref{1})). Applying to $f^{(0)}$ the constructions from the proof of Theorem~\ref{algebraic-rational} (cf. (\ref{3})) and from the proof of Theorem~\ref{tropical-rational}~(i), we complete the proof of (ii).

(iii). The proofs for both cases $f$ being a tropical increasing rational function or a tropical polynomial go similarly.

Let $f=g_1\circ \cdots \circ g_k$ where each $g_i,\, 1\le i \le k$ is a tropical increasing rational binomial (respectively, a tropical binomial) with slopes $b_i,\, c_i,\, 1\le i\le k$ (respectively, $b_i>c_i$). Denote by $r_i,\, 1\le i\le k$ the unique tropical root of $g_i$. Partition $\RR$ into intervals with the end-points $(g_{i+1}\circ \cdots \circ g_k)^{-1}(r_i),\, 1\le i\le k$. In case of tropical polynomials $f$ all these end-points are the tropical roots of $f$. In case of tropical increasing rational functions $f$ all $g_i$ for which $(g_{i+1}\circ \cdots \circ g_k)^{-1}(r_i)$ being not a tropical root of $f$, give a contribution into $g_1\circ \cdots \circ g_k$ by multiplying all the slopes of its edges on the intervals by the same integer, so w.l.o.g. one can assume that each $(g_{i+1}\circ \cdots \circ g_k)^{-1}(r_i),\, 1\le i\le k$ is a tropical root of $f$.

For $1\le j\le n$ take the set $I_j$ of $1\le i\le k$ such that $(g_{i+1}\circ \cdots \circ g_k)^{-1}(r_i)$ is $j$-th root $t_j$ of $f$. Then $a_j/a_{j-1}=\prod _{i\in I_j} (c_i/b_i)$. Therefore, $q_j| \prod _{i\in I_j} b_i$. Since $\prod _i b_i =a_0$, we conclude that $(q_1\cdots q_n)|a_0$.

Conversely, let $(q_1\cdots q_n)|a_0$. Put integers $b_j:=q_j,\, 1\le j\le n-1,\, b_n:=a_0/(b_1\cdots b_{n-1})$ and $c_j:=b_ja_j/a_{j-1},\, 1\le j\le n$. 

Construct $g_n,\dots,g_1$ recursively. As a base of recursion take $g_n$ such that its unique tropical root coincides with $t_n$ (observe that $g_n$ is defined uniquely up to an additive shift, in other words, one can replace $g_n$ by $g_n+e,\, e\in \RR$). Assume that $g_n,\dots,g_{m+1}$ are already constructed by recursion. Then take $g_m$ such that its unique tropical root equals $(g_{m+1}\circ \cdots \circ g_n)(t_m)$. At the very last step of recursion we adjust $g_1$ by a suitable additive shift to make $g_1\circ \cdots \circ g_n$ coincide with $f$ at one (arbitrary) point. Hence $f= g_1\circ \cdots \circ g_n$. $\Box$

%(iii). Now let a tropical increasing rational trinomial $f$ have slopes $a,\, b,\, c>0$ of its edges. Assume that $f=g\circ h$ for some tropical rational binomials (clearly, this is the only possible shape of non-trivial decompositions of $f$). Obviously, $g,\, h$ are increasing. Denote by $x_0<x_1$ the tropical roots of $f$ (then the slope of $f$ on the interval $(-\infty,\, x_0]$ equals $a$, the slope on the interval $[x_0,\, x_1]$ equals $b$, and the slope on the interval $[x_1,\, \infty)$ equals $c$. One of the roots $x_0,\, x_1$ is the unique tropical root of $h$ (let it be $x_1$ for definiteness). Thus, $h$ coincides with a linear function $rx+r_0$ on the interval $(-\infty,\, x_1]$ and with a linear function $qx+q_0$ on the interval $[x_1,\, \infty)$ for suitable integers $r,\, q>0$. Therefore, $g$ coincides with a linear function $sx+s_0$ on the interval $(-\infty,\, rx_0+r_0]$ and with a linear function $tx+t_0$ on the interval $[rx_0+r_0,\, \infty)$. Hence $a=sr,\, b=tr,\, c=tq$, i.~e. $b|(a\cdot c)$.

%Conversely, if $b|(a\cdot c)$ then put $r:=GCD(a,\, b),\, t:=b/r$. We have $t|c$ and put $s:=a/r,\, q:=c/t$. One can take $r_0:=0,\, q_0:=(r-q)x_1,\, s_0:=f(x_0)-ax_0,\, t_0:=-bx_0+f(x_0)$. $\Box$ 

\begin{remark}
The algorithms designed in sections~\ref{two}, \ref{three} have polynomial complexity since after each procedure yielding a composant (cf. (\ref{1}), (\ref{2}), (\ref{3})) either the number of blocks of edges or the number of edges drops at least by one.
\end{remark}

\section{Tropical polynomials with commuting composition}\label{four}

Let $f,\, g$ be tropical polynomials without free terms (some statements below hold also for more general tropical increasing algebraic rational functions).
% (in particular, one can consider tropical polynomials without free terms). 
In this section we give a criterion when $f\circ g=g\circ f$. Note that the inverse $f^{-1}$  (i.~e. $f\circ f^{-1}=Id$ equals the identity function) is  a tropical increasing algebraic rational function. Denote by $f^k:= f\circ \cdots \circ f$ the $k$ times iteration of $f$. We agree that $f^0:=Id$. Note that tropical increasing algebraic rational functions constitute a group with respect to the composition.

\begin{remark}\label{fixed-point}
Let $f,\, g$ be tropical increasing algebraic rational functions and $f\circ g=g\circ f$ hold. We call $x$ a {\it fixed point} of $f$ if $f(x)=x$. The set $F_f\subset \RR$ of fixed points is a finite union of disjoint closed intervals $\{[x_i,\, y_i]\}_i$ (including isolated points, i.~e. $x_i=y_i$). Since $f\circ g(x)=g\circ f(x)=g(x)$ we conclude that $g(F_f)=F_f$, therefore $g(x_i)=x_i,\, g(y_i)=y_i$ for all $i$ since $g$ is increasing. 

Observe that either $g(x)=x$ for any point $y_i<x<x_{i+1}$, either $g(x)<x$ for any point $y_i<x<x_{i+1}$ or $g(x)>x$ for any point $y_i<x<x_{i+1}$. Indeed, otherwise consider the set of fixed points $F_g\cap [y_i,\, x_{i+1}]$, and arguing as above in the previous paragraph we get $f(F_g\cap [y_i,\, x_{i+1}])=F_g\cap [y_i,\, x_{i+1}]$, again $f(x)=x$ for any end-point of an interval of $F_g\cap [y_i,\, x_{i+1}]$, which contradicts the choice of $y_i,\, x_{i+1}$, unless $F_g\supset (y_i,\, x_{i+1})$, in other words $g(x)=x$ for any point $y_i<x<x_{i+1}$. We allow intervals with $\pm \infty$ end-points.  

In the case of tropical polynomials $f,\, g$ without free terms there is at most one end-point of the intervals of fixed points of $f,\, g$, which we denote by $x_0$, due to the convexity of $f,\, g$ and taking into the account that the slopes of edges of $f,\, g$ are greater or equal than 1. Thus, there are at most two intervals $(-\infty,\, x_0],\, [x_0,\, \infty)$ or just one interval $(-\infty,\, \infty)$ when $x_0=\infty$.

\end{remark}

\begin{theorem}\label{commuting}
Tropical polynomials $f,\, g$ commute: $f\circ g=g\circ f$ iff 
either $x_0=\infty$ and $f=x+c_1,\, g=x+c_2, x\in \RR$ for some $c_1,\, c_2\in \RR$ or $-\infty<x_0<\infty$ and the following is valid.

There exists a tropical increasing algebraic rational function $h$ such that $h(x_0)=x_0$  (see Remark~\ref{fixed-point}) and 

$\bullet$ either $f=h^{p},\, g=h^{q}$ for suitable non-negative integers $p,\, q$

$\bullet$ or $f=ax+x_0(1-a),\, g=bx+x_0(1-b)$ for suitable integers $a,\, b\ge 1$

\noindent holds on the interval $[x_0,\, \infty)$.  Similarly,

$\bullet$ either $f=h^{k},\, g=h^{m}$ for suitable non-negative integers $k,\, m$

$\bullet$ or $f=dx+x_0(1-d),\, g=ex+x_0(1-e)$ for suitable integers $d,\, e\ge 1$

\noindent holds on the interval $(-\infty,\, x_{0}]$,

%either $f=h^{p},\, g=h^{q}$ on the interval $(x_0,\, \infty)$

%\noindent for suitable non-negative integers $k,\, m,\, p,\, q$, unless  $f=x+c_1,\,% g=x+c_2$ for some $c_1,\, c_2\in \RR$.
\end{theorem}

\begin{remark}
Note that $h$ is not necessary a tropical polynomial.
\end{remark}

 {\bf Proof}. In one direction, namely when either such appropriate $h$ does exist or $c_1,\, c_2$ do exist, obviously $f\circ g=g\circ f$ holds.

From now on let $f\circ g=g\circ f$. If $f(x)=x$ (respectively, $g(x)=x$) for any $x\ge x_0$ one can put $h:=g,\, f=h^0$ (respectively, $h:= f,\, g=h^0$) on the interval $(x_0,\, \infty)$. Thus, from now on we suppose that $f(x)>x,\, g(x)>x$ for any $x>x_0$ (cf. Remark~\ref{fixed-point}). We construct (increasing) $h$ on the interval $(-\infty,\, x_0)$ and separately on the interval $(x_0,\, \infty)$ such that $h(x_0)=x_0$ and after that glue them together and obtain a tropical increasing algebraic rational function $h$ required in Theorem~\ref{commuting}.

\begin{lemma}\label{value}
Let $f,\, g$  be tropical increasing algebraic rational functions, $f\circ g=g\circ f$ and for some point $y_i<x<x_{i+1}$ it holds $f(x)=g(x)$. Then $f$ coincides with $g$ on the interval $(y_i,\, x_{i+1})$.
\end{lemma} 

{\bf Proof of Lemma~\ref{value}}. Since neither $f$ nor $g$ has a fixed point in the interval $(y_i,\, x_{i+1})$  one can assume for definiteness that $f(y)>y,\, g(y)>y$ for any $y_i<y<x_{i+1}$ (see Remark~\ref{fixed-point}). For each integer $k$ we have $f(f^k(x))=g(f^k(x))$. The increasing sequence $x<f(x)<f^2(x)<\cdots$ tends to $x_{i+1}$ taking into the account that $F_f\cap (y_i\, x_{i+1})=\emptyset$. Therefore, the right-most edges of $f$ and $g$ on the interval $(y_i,\, x_{i+1})$ coincide. Suppose that $f$ and $g$ do not coincide on the interval $(y_i,\, x_{i+1})$.

Take the left-most point $z_0\in (y_i,\, x_{i+1})$ such that $f(y)=g(y)$ for any $z_0\le y<x_{i+1}$. Hence on a sufficiently small interval $[z,\, z_0]$ function $f$ (respectively, $g$) is linear $ax-az_0+f(z_0)$ (respectively, $bx-bz_0+f(z_0)$) and $a\neq b$. Since $x_{i+1}>f(z_0)=g(z_0)>z_0$ and due to the choice of $z_0$ there exists a linear function $cx+d$ such that on a sufficiently small interval $[z,\, z_0]$ the composition $f\circ g$ coincides with the linear function $cbx-cbz_0+cf(z_0)+d$, while the composition $g\circ f$ on $[z,z_0]$ coincides with the linear function $cax-caz_0+cf(z_0)+d$, which contradicts to the commutativity $f\circ g=g\circ f$ and proves Lemma~\ref{value}. $\Box$ \vspace{2mm}

Fix an interval $(y_i,\, x_{i+1})$ for the time being and denote by $T_f$ the set of tropical roots of $f$. We considered the case when $f(x)=x$ for any $x\in (y_i,\, x_{i+1})$ or when $g(x)=x$ for any $x\in (y_i,\, x_{i+1})$ above, so we assume that either $f(x)>x$ for any $x\in (y_i,\, x_{i+1})$ or $f(x)<x$ for any $x\in (y_i,\, x_{i+1})$, and  either $g(x)>x$ for any $x\in (y_i,\, x_{i+1})$ or $g(x)<x$ for any $x\in (y_i,\, x_{i+1})$ (cf. Remark~\ref{fixed-point}). First, we study the case $(T_f\cup T_g)\cap (y_i,\, x_{i+1}) =\emptyset$. Since $g(y_i)=f(y_i)=y_i,\, g(x_{i+1})=f(x_{i+1})=x_{i+1}$ we conclude that one or both end-points of the interval $(y_i,\, x_{x+1})$ equal $\pm \infty$. 

When both $y_i=-\infty,\, x_{i+1}=\infty$, we have $f=x+c_1,\, g=x+c_2$ for some $c_1,\, c_2 \in \RR$.

If $y_i\in\RR,\, x_{i+1}=\infty$ (the case $y_i=-\infty,\, x_{i+1}\in \RR$ is analyzed in a similar way) then $f$ (respectively, $g$) coincides on the interval $(y_i,\, \infty)$ with a linear function $ax-ay_i+y_i$ (respectively, $bx-by_i+y_i$) for suitable rationals $a,\, b>0$ which establishes Theorem~\ref{commuting} in the case $(T_f\cup T_g)\cap (y_i,\, x_{i+1}) =\emptyset$.

%$a=p/l,\, b=q/l$ with positive integers $p,\, q,\, l$. Put $h:=x/l-y_i/l+y_i$, and we get $f=h^p,\, g=h^q$ on the interval $(y_i,\, \infty)$, which proves Theorem~\ref{commuting} in the case of the interval $(y_i,\, \infty)$. 

From now on we again assume  $f,\, g$ to be tropical polynomials and let $(T_f\cup T_g) \cap (-\infty,\, x_0) \neq \emptyset$ (cf. Remark~\ref{fixed-point}).

\begin{lemma}\label{group}
(i) If $h_1,\, h_2$ are tropical  algebraic rational functions and $x\in T_{h_1\circ h_2}$ then either $x\in T_{h_2}$ or $h_2(x)\in T_{h_1}$. For tropical polynomials $f,\, g$ the converse is true: if either $x\in T_g$ or $g(x)\in T_f$ then $x\in T_{f\circ g}$;

(ii) let $f\circ g=g\circ f$. If $x\in T_g$  then either $f^{-1}(x)\in T_g$ or $g\circ f^{-1} (x)\in T_f$;

(iii) let $f\circ g=g\circ f$. If $x\in T_f\setminus T_g$ then $g(x)\in T_f$.
\end{lemma} 

{\bf Proof of Lemma~\ref{group}}. (i). For the converse statement the convexity of $f,\, g$ is used.

(ii). Denote $y:= f^{-1}(x)$. Due to (i)  $y\in T_{g\circ f}$, hence either $y\in T_g$ or $g(y)\in T_f$ again due to (i) and taking into the account that $f\circ g=g\circ f$;

(iii) Since due to (i) $x\in T_{g\circ f} =T_{f\circ g}$ we conclude that $g(x)\in T_f$ again by means of (i). $\Box$

\vspace{2mm}

Consider a directed graph $G$ with the nodes being the points from $(T_f\cup T_g) \cap (-\infty,\, x_0)$ and the arrows according to Lemma~\ref{group} as follows (recall that $G$ is not empty, the case of empty $G$ was studied above). From every node $x\in T_g \cap (-\infty,\, x_0)$ there is an arrow labeled by $f^{-1}$ to the node $f^{-1}(x)$, provided that $f^{-1}(x) \in T_g$, and there is an arrow labeled by $g\circ f^{-1}$ to the node $g\circ f^{-1}(x)$, provided that $g\circ f^{-1}(x) \in T_f$ (observe that $f^{-1}(x),\, g\circ f^{-1}(x) \in (-\infty,\, x_0)$). In addition, there is an arrow labeled by $g$ from every node $x\in T_f\setminus T_g$ to the node $g(x)\in T_f$ (again $g(x) \in (-\infty,\, x_0)$).

There is a cycle in $G$ (due to Lemma~\ref{group}), let it contain a node $x$. Denote by $t$ the composition of the labels of the arrows (starting with $x$) in this cycle. Then $t(x)=x$ and one can represent $t=g^s\circ f^{-r}$ (taking into the account that $f\circ g=g\circ f$) for some non-negative integers $s,\, r$ at least one of which being positive. Observe that in fact, $s,\, r>0$ since $f(x_1)<x_1,\, g(x_1)<x_1$ for any $x_1<x_0$ (cf. Remark~\ref{fixed-point}).

Hence $g^s(x)=f^r(x)$. Lemma~\ref{value} implies that $g^s$ coincides with $f^r$ on the interval $(-\infty,\, x_0)$. Denote $n:=GCD(s,\, r)$, then $(g^{s/n}\circ f^{-r/n})^n=Id$ on the interval $(-\infty,\, x_0)$. The function $u:=  g^{s/n}\circ f^{-r/n}$ is increasing piece-wise linear. Therefore, one can partition $\RR$ into a finite number of intervals (including unbounded ones) such that on each of these intervals $[y_0,\, y_1]$ it holds $u(y_0)=y_0,\, u(y_1)=y_1$, and either $u(y)>y$ for any $y_0<y<y_1$, either $u(y)<y$ for any $y_0<y<y_1$ or $u(y)=y$ for any $y_0<y<y_1$ (cf. Remark~\ref{fixed-point}). Hence $u=Id$, i.~e. $g^{s/n}=f^{r/n}$ on the interval $(-\infty,\, x_0)$.

For appropriate positive integers $i,\, j$ it holds $1=-i(s/n)+j(r/n)$. Consider a tropical increasing algebraic rational function $h:=g^j\circ f^{-i}$. Then

$h^{r/n}=g^{jr/n}\circ f^{-ir/n}=g^{jr/n}\circ g^{-is/n}=g$;

$h^{s/n}=g^{js/n}\circ f^{-is/n}=f^{jr/n}\circ f^{-is/n}=f$

\noindent on the interval $(-\infty,\, x_0)$. 

In a similar way one produces $h$ on the interval $(x_0,\, \infty)$, provided that $x_0<\infty$. This completes the proof of Theorem~\ref{commuting}. $\Box$

\begin{remark}
It would be interesting to give a criterion for commuting tropical increasing algebraic rational functions, and more generally, for tropical non-monotone algebraic rational functions.
\end{remark}

\begin{example}
We exhibit an example of a commuting pair of increasing tropical rational functions $f,\, g$ (defined on the interval $[x_0,\, \infty),\, f(0)=g(0)=0$, thereby $x_0=0$, see Theorem~\ref{commuting}) not satisfying the conclusion of Theorem~\ref{commuting}.

Pick a real $0<t$, integers $a>\alpha>1,\, b\ge 1$ such that $a\neq b,\, a|(b\alpha)$. denote $u:=\alpha t,\, v:=at,\, w:=\alpha a t$. Define $g$ to be piece-wise linear whose graph on $[x_0,\, \infty)$ consists of three edges having the slopes $a,\, b,\, a$, respectively, and with the tropical roots at the points $u,\, w$ (this defines $g$ uniquely). Similarly, define $f$ whose graph also has three edges with the slopes $\alpha,\, (b\alpha)/a,\, \alpha$, respectively, and with the tropical roots $v,\, w$.
Then the graph on $[x_0,\, \infty)$ of the increasing tropical rational function $f\circ g=g\circ f$ has three edges as well with the slopes $\alpha a,\, \alpha b,\, \alpha a$, respectively, and with tropical roots $t,\, w$. 

In case when $f^k=g^m$ (cf. Theorem~\ref{commuting}) we get that $\alpha^k=a^m$. Thus, if there are no such integers $k,\, m$ the conclusion of Theorem~\ref{commuting} for $f,\, g$ is not fulfilled. 
\end{example}

\section{Tropical polynomial and rational  parametrizations}\label{five}

We call a {\it polygonal line} $L\subset \RR^n$ with $k+1$ intervals a sequence of intervals with endpoints $v_1,\dots,v_k \in \QQ^n$ such that $i$-th interval has endpoints $v_i,\, v_{i+1}$ for $1\le i\le k-1$, while unbounded $0$-th interval (a ray) has $v_1$ as its right endpoint, and unbounded $k$-th interval (a ray) has $v_k$ as its left endpoint. The {\it vector of slopes} of $i$-th interval, $1\le i\le k-1$ is defined as $(a_{i,1},\dots,a_{i,n}):=v_{i+1}-v_i$, similarly one can define a vector of slopes $(a_{0,1},\dots,a_{0,n})$ of $0$-th and $(a_{k,1},\dots,a_{k,n})$ of $k$-th intervals, respectively (we assume that the latter two vectors of slopes are also rational). Note that the vector of slopes is determined up to a positive factor.

We say that tropical rational functions $f_1,\dots,f_n$ in one variable $t$ provide a {\it tropical rational parametrisation} of $L$ if the map $(f_1,\dots,f_n):\RR \to L$ is a bijection and (for definiteness) $(f_1,\dots,f_n)^{-1}(v_i)< (f_1,\dots,f_n)^{-1}(v_{i+1}),\, 1\le i\le k-1$. In particular, $\{v_1,\dots,v_k\}$ coincides with the set of all the tropical roots of  $f_1,\dots,f_n$, i.~e. the points where one of the functions $f_1,\dots,f_n$ is not smooth. We suppose w.l.o.g. that one can't discard any $v_i,\, 1\le i\le k$ while keeping the propery of $L$ to be a polygonal line. When $f_1,\dots,f_n$ are tropical polynomials (respectively, Laurent polynomials), we talk about {\it tropical polynomial} (respectively, {\it Laurent polynomial)   parametrisation} of $L$.

In case if $L$ is a subset of a tropical curve one can treat a parametrisation of $L$ as a parametrisation of the tropical curve (cf. \cite{G}) since a parametrisation provides a parametric family of solutions of a  system of tropical equations.

\begin{example}
Let $T\subset \RR^2$ be a tropical curve (a tropical line) defined by a tropical polynomial $\min\{x,\, y,\, 0\}$. Then $L\subset T$ consisting of two rays $\{x=0\le y\} \cup \{y=0\le x\}$ admits a tropical rational parametrisation with
$f_1:=-\min\{t,\, 0\},\, f_2:=-\min\{-t,\, 0\}$.
\end{example} 

\begin{proposition}\label{parametrisation}
A polygonal line $L$ has 

(i) always a tropical rational parametrisation; \vspace{2mm}

(ii) a tropical polynomial parametrisation iff $a_{i,j}\ge 0,\, 0\le i\le k,\, 1\le j\le n$, and  $a_{i,j}=0$ implies $a_{l,j}=0$ for all $l\ge i$; \vspace{2mm}

(iii) a tropical Laurent polynomial parametrisation iff

\noindent $\bullet$ $a_{i,j}<0$ implies $a_{i+1,j}<0$;

\noindent $\bullet$ $a_{i+1,j}>0$ implies $a_{i,j}>0$;

\noindent $\bullet$ $a_{i,j_0}>0,\, a_{i+1,j_0}>0,\, a_{i,j}<0,\, a_{i+1,j}<0$ imply
$a_{i,j_0}/a_{i,j}\le a_{i+1,j_0}/a_{i+1,j}$ 

\noindent for all $0\le i\le k-1,\, 1\le j\neq j_0 \le n$.   
\end{proposition}

{\bf Proof}. (i) We have to construct tropical univariate rational functions $f_1,\dots,f_n$. First we construct piece-wise linear functions $g_1,\dots,g_n$ with rational slopes (in \cite{G} such functions are called tropical Newton-Puiseux rational functions). As a set of tropical roots of $g_1,\dots,g_n$ we take points $1,\dots,k$. The vector of the values of $g_1,\dots,g_n$ at point $i$ we put $v_i,\, 1\le i\le k$. Thereby, $g_1,\dots,g_n$ are defined on interval $[1,\, k]\subset \RR$. To extend $g_1,\dots,g_n$ to interval $(-\infty,\, 1]$ (respectively, $[k,\, \infty)$) use the vector of the slopes of $0$-th (respectively, $k$-th) interval of $L$.

To proceed to tropical rational functions $f_1,\dots,f_n$ (so, piece-wise linear functions with integer slopes), denote by $M$ the least common multiple of all the denominators of the slopes of $g_1,\dots,g_n$ (i.~e. the slopes of $L$). As the set of tropical roots of $f_1,\dots,f_n$ take points $1/M,\dots,k/M$. The vector of the values at point $i/M$ we put $v_i,\, 1\le i\le k$. In other words, the corresponding slopes of $f_1,\dots,f_n$ are obtained from the corresponding slopes of $g_1,\dots,g_n$ multiplying by $M$. Satisfying also the latter condition, one extends $f_1,\dots,f_n$ to intervals $(-\infty,\, 1/M]$ and $[k/M,\, \infty)$. \vspace{1mm}

(ii) If $f_1,\dots,f_n$ constitute a tropical polynomial parametrization of $L$ then since the slopes of each $f_j,\, 1\le j\le n$ (being a convex function) are non-increasing non-negative integers we get the conditions stated in (ii).

Conversely, if the latter conditions are fulfilled one can recursively on $i$ choose positive rationals $c_0=1,\, c_1,\dots,c_k$ in such a way that $c_{i+1}\cdot a_{i+1,j}\le c_i \cdot a_{i,j},\, 0\le i\le k-1,\, 1\le j\le n$ taking each $c_{i+1}$ to be the maximal possible among satisfying the latter inequalities. Therefore, one can take $c_i\cdot a_{i,j},\, 1\le j\le n$ as the slopes of (Newton-Puiseux polynomials \cite{G}, i.~e. convex piece-wise linear functions with rational non-negative slopes) $g_j,\, 1\le j\le n$ with the tropical roots at points $1,\dots,k$. Then as at the end of the proof of (i) one can obtain tropical polynomials $f_j,\, 1\le j\le n$ with the non-negative integer slopes $M\cdot c_i\cdot a_{i,j},\, 0\le i\le k,\, 1\le j\le n$ and with the tropical roots at points $1/M,\dots,k/M$. Then $f_1,\dots,f_n$ provide a required parametrization of $L$. \vspace{1mm}

(iii) If there exists a tropical Laurent polynomial parametrization $f_1,\dots,f_n$ of $L$ then the slopes $b_{i,j},\, 0\le i\le k,\, 1\le j\le n$ of $f_j,\, 1\le j\le n$, respectively, being integers fulfil the conditions $b_{i,j}\ge b_{i+1,j},\, 0\le i\le k-1,\, 1\le j\le n$. On the other hand, there exist positive rationals $c_0,\dots,c_k$ such that 
$b_{i,j}=c_i\cdot a_{i,j},\, 0\le i\le k,\, 1\le j\le n$. This entails the conditions from (iii).

Conversely, let the conditions from (iii) be fulfilled. Construct positive rationals $c_0=1,\, c_1,\dots, c_k$ such that $c_i\cdot a_{i,j}\ge c_{i+1}\cdot a_{i+1.j},\, 0\le i\le k-1,\, 1\le j\le n$ by recursion on $i$. Assume that $c_0=1,\, c_1,\dots,c_i$ are already constructed. Take the maximal possible $c_{i+1}>0$ such that $c_i\cdot a_{i,j}\ge c_{i+1}\cdot a_{i+1.j}$ for all $1\le j\le n$ such that $a_{i,j}>0,\, a_{i+1,j}>0$. Then for suitable $j_0$ for which $a_{i,j_0}>0,\, a_{i+1,j_0}>0$ it holds $c_i\cdot a_{i,j_0}=c_{i+1}\cdot a_{i+1,j_0}$. For every $1\le j\le n$ for which $a_{i,j}<0,\, a_{i+1,j}<0$ the condition from (iii)   $a_{i,j_0}/a_{i,j}\le a_{i+1,j_0}/a_{i+1,j}$ implies $c_{i+1}\cdot a_{i+1,j} \le c_i\cdot a_{i,j}$.

Thus, as in (i), (ii) one first constructs piece-wise linear functions $g_j,\, 1\le j\le n$ with rational non-increasing slopes $c_i\cdot a_{i,j},\, 1\le j\le n$ and with the tropical roots at points $1,\dots,k$. Denote by $M$ the common denominator of these slopes and construct tropical Laurent polynomials $f_1,\dots,f_n$ with the slopes obtained from the slopes of $g_j,\, 1\le j\le n$ multiplying them by $M$ and with the tropical roots $1/M,\dots,k/M$. Then $f_1,\dots,f_n$ provide a required parametrization of $L$. $\Box$

\begin{remark}
One can construct the required parametrizations in Proposition~\ref{parametrisation} within polynomial complexity following the proofs of (i), (ii), (iii).
\end{remark}

It would be interesting to extend parametrizations from 1-dimensional polygonal lines to multidimensional polyhedral complexes.
\vspace{2mm}

{\bf Acknowledgements}. The author is grateful to the grant RSF 16-11-10075
and to MCCME for inspiring atmosphere.


\begin{thebibliography}{99}
%\bibitem{G}
%D.~Grigoriev.
%\newblock On a tropical dual Nullstellensatz.
%\newblock {\em Adv. Appl. Math.}, 48:457--464, 2012.

\bibitem{Gathen}
J.~von zur Gathen.
\newblock Functional decomposition of polynomials: the wild case.
\newblock {\em  J. Symbolic Comput.}, 10:437--452, 1990.

\bibitem{Rubio}
J.~von zur Gathen, J.~Gutierrez and R.~Rubio.
\newblock Multivariate polynomial decomposition.
\newblock {\em  Appl. Algebra Engrg. Comm. Comput.}, 14:11--31, 2003. 

\bibitem{G}
D.~Grigoriev.
\newblock Tropical Newton-Puiseux polynomials.
\newblock {\em Lect. Notes Comput. Sci.}, 11077:177--186, 2018.

\bibitem{GP}
D.~Grigoriev and V.~Podolskii.
\newblock Tropical combinatorial Nullstellensatz and fewnomials testing.
\newblock {\em Lect. Notes Comput. Sci.}, 10472:284--297, 2017.

\bibitem{Kozen}
D.~Kozen, S.~Landau and R.~Zippel.
\newblock Decomposition of algebraic functions.
\newblock {\em  J. Symbolic Comput.}, 22:235--246, 1996.

\bibitem{MS}
D.~Maclagan and B.~Sturmfels.
\newblock {\em Introduction to Tropical Geometry:}, volume 161 of {\em Graduate
  Studies in Mathematics}.
\newblock American Mathematical Society, 2015.

\bibitem{neural}
G.~F.~Mont\'ufar, R.~Pascanu, K.~Cho and Y.~Bengio.
\newblock On the number of linear regions of deep neural networks.
\newblock {\em Advances in Neural Information Processing Systems 27}, Montreal,
2924--2932, 2014.

\bibitem{P}
F.~Pakovich.
\newblock Semiconjugate Rational Functions: A Dynamical Approach.
\newblock {\em Arnold Math. J.}, 4:59--68, 2018.

\bibitem{R22}
J.~Ritt.
\newblock Prime and composite polynomials.
\newblock {\em  Trans. Amer. Math. Soc.}, 23:51--66, 1922.

\bibitem{R23}
J.~Ritt.
\newblock Permutable rational functions.
\newblock {\em  Trans. Amer. Math. Soc.}, 25:399--448, 1923.

\end{thebibliography}
\end{document}